\documentclass{amsart}
\usepackage{amssymb}
\usepackage{amsmath}

\newtheorem{theorem}{Theorem}
\newtheorem{Main theorem}{Main theorem}
\newtheorem{prop}{Proposition}
\newtheorem{lemma}{Lemma}
\newtheorem{corollary}{Corollary}

\newtheorem{conjecture}{Conjecture}
\newtheorem{example}{Example}
\newcommand{\N}{\ensuremath{\mathbb{N}}}
\newcommand{\Z}{\ensuremath{\mathbb{Z}}}
\newcommand{\Q}{\ensuremath{\mathbb{Q}}}
\newcommand{\R}{\ensuremath{\mathbb{R}}}
\newcommand{\C}{\ensuremath{\mathbb{C}}}
\newcommand{\Pro}{\ensuremath{\mathbb{P}}}
\begin{document}
\title{Oscillating integrals and Newton polyhedra}
 \author[J. Denef]{Jan Denef}
\address{Department of Mathematics\\
Katholieke Universiteit Leuven\\ Celestijnenlaan 200B\\ B-3001
Leuven\\ Belgium}
 \email{jan.denef@wis.kuleuven.ac.be}
 \urladdr{http://www.wis.kuleuven.ac.be/algebra/denef/}
\author[J. Nicaise]{Johannes Nicaise$^\dag$}
\address{Department of Mathematics\\
Katholieke Universiteit Leuven\\ Celestijnenlaan 200B\\ B-3001
Leuven\\ Belgium}
 \email{johannes.nicaise@wis.kuleuven.ac.be}
\thanks{$\dag$Research Assistant of the Fund for Scientific Research --
 Flanders (Belgium)(F.W.O.)}
\author[P. Sargos]{Patrick Sargos}
\address{Institut Elie Cartan\\ Universit\'e Henri Poincar\'e,
Nancy 1\\B.P. 239\\ 54506 Vandoeuvre-les-Nancy\\ France}
\email{sargos@iecn.u-nancy.fr}
 \maketitle
\section{Introduction}
Oscillating integrals are integrals of the form
$$\int_{\R^{n}}e^{itf(x)}\varphi(x)dx\,.$$
They frequently occur in applied mathematics and mathematical
physics. In this article, we investigate their asymptotic
behaviour when the parameter $t$ tends to infinity, in terms of
the geometry of the Newton polyhedron of the phase $f$. It is
well-known that the greatest contributions to this asymptotic
behaviour arise from the critical points of $f$: if $f$ is regular
on the support of $\varphi$, subsequent oscillations will more or
less cancellate each other as $t$ grows bigger and the integrand
starts oscillating faster, so that the integral tends to zero more
rapidly than any power of $t$. This phenomenon is called the
principle of the stationary phase. When $f$ has only
non-degenerate critical points, we can apply Morse's lemma to give
a description of the asymptotic expansion of the integral; see
\cite{Arnold}. In the present paper, we will consider a much
larger class of phase functions $f$: real analytic functions which
are non-degenerate with respect to their Newton polyhedron.

  The key result is Theorem 1 in section 5, which yields,
together with formula (1) in section 3, an expression of
$\mu(\varphi)$ in terms of principal value integrals, where
$\mu(\varphi)$ is the coefficient of the expected leading term in
the asymptotic expansion of our oscillating integral. A similar -
but more complicated - expression was given in \cite{DeSaLa} in terms
of a different kind of principal value integrals. A direct
consequence of Theorem 1 is Corollary 1, which states that the
coefficients $\mu(\varphi)$ for $f$ and $f_{\tau_{0}}$ differ only
by an easy nonzero factor. The much simpler polynomial
$f_{\tau_{0}}$, as defined in the next section, is obtained by
omitting all monomials of $f$ whose exponents do not lie on the
face $\tau_{0}$. Here $\tau_{0}$ is the smallest face of the
Newton polyhedron of $f$ intersecting the diagonal.

  As a first application, we give in Section \ref{conj} a
more transparent proof of the fact that $\mu(\varphi)=0$ whenever
$\tau_{0}$ is unstable. This result, conjectured by Denef and
Sargos, was first proven in \cite{DeSaLa}. As a second
application, we give a very explicit formula for $\mu(\varphi)$ in
terms of gamma functions, assuming that $\tau_{0}$ is a simplex of
codimension $1$, whose vertices are the only integral points on
$\tau_{0}$ corresponding to monomials of $f$. This is done in
Section \ref{explicit}.

In Section \ref{sectioncomplex}, we develop an analogous residue
formula for the complex local zeta function. This allows us to
give, in section \ref{sectionstability}, a partial proof of the
stability conjecture of Denef and Sargos, using a theorem of
Loeser on the relation between the spectrum of a complex
polynomial, and the poles of the associated complex local zeta
function.

  We conclude this introduction with an overview of the structure of the paper.
Section 2 contains some preliminaries concerning Newton polyhedra,
while section 3 deals with some known results about the asymptotic
behaviour of oscillating integrals. In section 4, we summarize the
analytic construction of toric varieties, which will be used to
desingularize $f$ at the origin. Section 5 establishes the real
residue formula. In section 6, we prove that $\mu(\varphi)$
vanishes if $\tau_{0}$ is unstable. The explicit formula for
$\mu(\varphi)$ is deduced in section \ref{explicit}. We give a
complex residue formula in Section \ref{sectioncomplex}, while in
section \ref{sectionstability}, we discuss the stability
conjecture in the complex case.

 We use the standard notation
$\R_{+}$ for the set of positive real numbers,
 and $\R_{0}$, $\C_0$, resp. $\N_0$, for
$\R\setminus \{0\}$, $\C\setminus\{0\}$, resp. $\N\setminus\{0\}$.
 To avoid confusion: $\R_{0}^{n}$ means
$(\R_{0})^{n}$, and $\C_{0}^{n}$ means
$(\C_{0})^{n}$.

\section{Newton polyhedra: some terminology}\label{poly}
Let $f$ be a nonconstant real analytical function in $n$ variables
$x=(x_{1},\ldots,x_{n})$ on an neighbourhood of $0\in \R^{n}$,
which has a critical point at the origin and satisfies $f(0)=0$.
Let $f(x)=\sum_{\alpha\in\N^{n}}a_{\alpha}x^{\alpha}$ be the
Taylor series of $f$ about the origin. We define the support
supp($f$) of $f$ to be the set of exponents $\alpha\in\N^{n}$ for
which $a_{\alpha}\neq 0$. The Newton polyhedron $\Gamma(f)$ of $f$
is the convex hull of supp($f$)$+\R_{+}^{n}$. For each subset
$\gamma$ of $\R^{n}_{+}$ we note by $f_{\gamma}$ the function
$f_{\gamma}(x)=\sum_{\alpha\in
\gamma\cap\N^{n}}a_{\alpha}x^{\alpha}$. When $\Gamma_{c}$ is the
union of the compact faces of $\Gamma(f)$, we call
$f_{\Gamma_{c}}$ the principal part of $f$. Let
$(t_{0},\ldots,t_{0})$ be the intersection point of the union of
faces of $\Gamma(f)$ with the diagonal $x_{1}=\cdots=x_{n}$, and
$\tau_{0}$ the smallest face of $\Gamma(f)$ containing this point.
Then $s_{0}$ denotes the value $-1/t_{0}$ and $\rho$ denotes the
codimension of $\tau_{0}$ in $\R^{n}$. We will sometimes write
$\tau_0(f)$ and $s_0(f)$ instead of $\tau_0$ and $s_0$, to make
$f$ explicit.

  We say that $f$ is non degenerate over $\R$ with respect
to $\Gamma(f)$ if the following holds: for each compact face
$\gamma$ of $\Gamma(f)$ the polynomial $f_{\gamma}$ has no
critical points in $\R_{0}^{n}$. "Almost all" phase functions $f$
are non degenerate with respect to $\Gamma(f)$. To specify this
assertion a little further: given a fixed Newton polyhedron
$\Gamma$, the principal parts of functions $f$ that are non
degenerate with respect to $\Gamma=\Gamma(f)$ form a
semi-algebraic subset of the space of principal parts with Newton
polyhedron $\Gamma$, and its complement is everywhere dense. From
now on, we will always assume that $f$ is non degenerate with
respect to its Newton polyhedron.

  To this Newton polyhedron we can associate a fan of
rational cones subdividing $\R^{n}_{+}$. The trace function
$l_{\Gamma}$ maps a vector $a$ in $\R^{n}_{+}$ to the value
$\min_{k\in \Gamma}<a,k>$, where $<,>$ denotes the inner product.
We define the trace of $a$ to be the compact face
$$\tau_{a}=\{k\in \Gamma_{c}\,|\,<a,k>=l_{\Gamma}(a)\}.$$ In
this way each $k$-dimensional compact face $\tau$ corresponds to a
$(n-k)$-dimensional cone $\dot{\tau}$ consisting of al the vectors
with trace $\tau$. Geometrically, this is the cone spanned by the
normal vectors of the facets of $\Gamma(f)$ containing $\tau$. It
is clear that these cones form a fan $\dot{\Gamma}$. We say that a
fan is subordinate to $\Gamma(f)$ if each cone of the fan is
contained in a cone of $\dot{\Gamma}$.

  An important role in this article is fulfilled by the
notion of instability. A face $\tau$ of $\Gamma(f)$ is unstable
over $\R$ with respect to the variable $x_{j}$, if $\tau$ is
contained in the region $\{y\in \R^{n}_{+}\,|\,0\leq y_{j}\leq
1\}$, but not entirely in the hyperplane defined by $x_{j}=0$, and
if furthermore, each compact face $\gamma$ of $\Gamma(f)$ that is
contained in the hyperplane $x_{j}=1$, is subject to the condition
that $f_{\gamma}$ has no zero in $\R_{0}^{n}$.

These notions can also be defined in the complex case.
Let $g$ be a non-constant
 complex analytic function on a neighbourhood of the origin in
$\C^{n}$, which has a critical point at the origin and satisfies $g(0)=0$.
The Newton polyedron $\Gamma(g)$ of $g$ at $0$ is defined as in the real case.
We say that $g$ is non degenerate over $\C$ with respect
to $\Gamma(g)$ if the following holds: for each compact face
$\gamma$ of $\Gamma(g)$ the polynomial $g_{\gamma}$ has no
critical points in $\C_{0}^{n}$, where $g_{\gamma}$
is defined in the same way as before. Again, "almost all"
complex phase functions $g$
are non degenerate with respect to $\Gamma(g)$.

  We conclude this section with some additional notation:
given a vector $\xi_{i}$ we mean by $\nu_{i}$ the sum of its
coordinates, and $N_{i}$ is short for $l_{\Gamma}(\xi_{i})$. If
the vector these notations are relating to is not explicitly
indicated it should be clear from the context which vector is
meant.
\section{Oscillating integrals}\label{osc}
 Let $\varphi$ be a
$\mathcal{C}^{\infty}$-function on $\R^{n}$ with support in a
sufficiently small neighbourhood of 0.

  It is well-known that the oscillating integral
$$I(t)=\int_{\R^{n}}e^{itf(x)}\varphi(x)dx$$ has for $t\rightarrow
\infty$ an asymptotic expansion
$$\sum_{p}\sum_{i=0}^{n-1}a_{p,i}(\varphi)t^{p}(\ln\,t)^{i}\qquad (*)$$
where $p$ runs through a finite number of arithmetic progressions,
not depending on the amplitude $\varphi$, that consist of negative
rational numbers. Since our objective is to study the asymptotic
behaviour of $I$, our primary interest goes out to the largest $p$
occurring in this expansion. Let $S$ be the set of tuples $(p,i)$
such that for each neighbourhood of $0$ in $\R^{n}$ there exists a
$\mathcal{C}^{\infty}$-function $\varphi$ with support in this
neighbourhood for which $a_{p,i}(\varphi)\neq 0$. We define the
oscillating index $\beta$ of $f$ to be the maximum of values $p$
for which we can find an $i$ so that $(p,i)$ belongs to $S$; the
maximum of these $i$ is called the multiplicity $\kappa$ of
$\beta$. The index $\beta$ contains information about the nature
of the singularity $0$ of $f$.

  In this paper, we will derive information about this
asymptotic expansion from the geometry of the Newton polyhedron
$\Gamma(f)$ of the phase function $f$, always assuming that $f$ is
non-degenerate with respect to $\Gamma(f)$ . Using the notation
introduced in the previous section, the expansion (*) can be
written as
$$\mu(\varphi)t^{s_{0}}(\log t)^{\rho-1}+O(t^{s_{0}}(\log
t)^{\rho-2}).$$ This is the main result in Varchenko's paper [8].
Moreover, it is known (as one can verify in [1]) that
$\beta=s_{0}$ and $\kappa=\rho-1$ if at least one of the following
additional conditions is satisfied:
$$\sharp\left\{
\begin{array}{ll} a) & s_{0}>-1
\\b) & f(x)\geq 0\ \mathrm{for\ all}\ x\in \R^{n}
\\c)& \tau_{0}\mathrm{\ is\ compact,\ } s_{0}\mathrm{\  is\ not\ an\ odd\ integer,\
and\ }f_{\tau_{0}}\mathrm{\ does\ not\ vanish\ in\ }
\R_{0}^{n}.\end{array}\right.
$$

  On $D=\{s\in\C\,|\,\Re(s)>0\}$ we can define the
functions
$$Z_{\pm}(s)=\int_{\R^{n}}f_{\pm}(x)^{s}\varphi(x)\,dx\ .$$ Here
$f_{+}=\mathrm{max}(f,0)$ and $f_{-}=\mathrm{max}(-f,0)$. It is
well known that these functions allow a meromorphic continuation
to the whole of $\C$, which we denote again by $Z_{\pm}(s)$.
 The study of the asymptotic behaviour of $I(t)$ can
be reduced to an investigation of the poles of $Z_{\pm}(s)$: the
terms in the development (*) are related to the singular part of
the Laurent expansion about the poles of $Z_{\pm}(s)$ (cf.
\cite{Arnold}).

  For a negative integer pole of $Z_{\pm}(s)$, it occurs
that terms in the singular part of the Laurent expansion at this
pole do not correspond to anything in the asymptotic expansion of
$I(t)$. This is why while studying the distribution $f_{\pm}^{s}$
we will always assume that $s_{0}\notin \Z$; the other case is
dealt with in the last paragraph of this section. When
$s_{0}\notin \Z$, $Z_{\pm}(s)$ can have a pole of order at most
$\rho$ at $s=s_{0}$. We define $\mu_{+}(\varphi)$ and
$\mu_{-}(\varphi)$ by requiring that
$$Z_{\pm}(s)=\frac{\mu_{\pm}(\varphi)}{(s-s_{0})^{\rho}}+O(\frac{1}{(s-s_{0})^{\rho-1}})$$
for $s\rightarrow s_{0}$. Using material in \cite{Arnold}, one obtains
the following expression for $\mu(\varphi)$ in terms of
$\mu_{+}(\varphi)$ and $\mu_{-}(\varphi)$:
\begin{equation}
\mu(\varphi)=\frac{1}{(\rho-1)!}\Gamma(-s_{0})[\,\mu_{+}(\varphi)e^{-\frac{i\pi
s_{0}}{2}}+\mu_{-}(\varphi)e^{\frac{i\pi s_{0}}{2}}\,]\ .
\end{equation}
 Note that, since $s_{0}\notin \Z$, this equality implies
that $\mu(\varphi)=0$ iff $\mu_{+}(\varphi)=\mu_{-}(\varphi)=0$.
This formula enables us to reduce the study of $\mu(\varphi)$ to
the study of the candidate pole $s_{0}$ of $Z_{\pm}(s)$. By the
relationship
$$Z_{\pm}(s)=\sum_{\theta \in
\{-1,1\}^{n}}\int_{\R^{n}_{+}}f_{\pm}(\theta_{1}x_{1},\ldots,\theta_{n}x_{n})^{s}\varphi(\theta
x)\,dx\, ,$$ it even suffices to investigate the properties of
$$\int_{\R^{n}_{+}}f_{\pm}(x)^{s}\varphi(x)\,dx\ .$$

  The question remains what happens when $s_{0}\in \Z$.
However, we can reduce this problem to the case $s_{0}\notin \Z$
by introducing an additional variable $y$: if we define a function
$f^{*}$ on $\R^{n}$ by
$f^{*}(x_{1},\ldots,x_{n},y)=f(x_{1},\ldots,x_{n})+y^{2}$ and put
$\varphi^{*}(x_{1},\ldots,x_{n},y)=\varphi(x_{1},\ldots,x_{n})\psi(y)$,
where $\psi$ is a test function on $\R$ satisfying $\psi(0)=1$,
and we define $\tau_{0}^{*}$, $\rho^{*}$, $s_{0}^{*}$ and
$\mu^{*}(\varphi^{*})$ in the obvious way, we obtain the following
properties:
\begin{itemize}
\item $f^{*}$ is non-degenerate over $\R$ with respect to its
Newton polyhedron
\item $s_{0}^{*}=s_{0}-\frac{1}{2}$
\item $\rho^{*}=\rho$
\item $\tau^{*}_{0}$ is the convex hull of $\tau_{0}$ and the
point $(0,\ldots,0,2)$ in $\R^{n+1}$
\item $\tau^{*}_{0}$ is stable iff $\tau$ is stable.
\end{itemize}
Moreover, $I^{*}(t)=I(t)\int_{\R}e^{ity^{2}}\psi(y)\,dy$. Since
the asymptotic expansion of this last factor equals
$$e^{\frac{i\pi}{4}}\sqrt{\frac{\pi}{t}}+O(t^{-\frac{3}{2}})\, ,$$
we conclude that
$$\mu(\varphi)=\frac{e^{-\frac{i\pi}{4}}}{\sqrt{\pi}}\mu^{*}(\varphi^{*})\ .$$
\section{Toric varieties}
Toric varieties form an important topic in algebraic geometry
because the geometric properties of this large class of varieties
are related to the combinatorial properties of the fans used to
construct them \cite{F}. Here we will introduce the analytical counterpart
of this construction, using an atlas with monomial transition
functions, following the approach in \cite{Arnold}.

  Let $L\subset \R^{n}$ be a lattice, e.g. $L=\Z^{n}$. A
cone in $\R^{n}$ is called rational if it can be generated by
vectors in $L$. We say the cone is simplicial if it can be
generated by a free set of vectors in $L$, and simple if this set
can be extended to a basis of $L$. Starting from a fan $F$ of
$L$-simple cones in $\R^{n}$, we will construct a real analytic
manifold $X_{L,F}$; this is the toric manifold associated to $L$
and $F$.

  We number once and for all the $1$-dimensional cones in
our fan $F$; this will allow us to speak of an ordered basis of
generators of a $n$-dimensional cone. The analytical structure of
$X_{L,F}$ is defined by giving an atlas for $X_{L,F}$, or more
specifically by giving a number of copies of $\R^{n}$ and the
transition functions between the parts of these copies that will
overlap once we identify these copies with open parts of
$X_{L,F}$.

  The charts $U_{\tau}$ in our atlas correspond to the
$n$-dimensional cones $\tau$ in our fan $F$, and an ordered basis
of generators of this cone provides standard coordinates on the
corresponding chart. Now we explain how you can travel from one
chart to another. Given two charts $U_{\tau_{1}}$ and
$U_{\tau_{2}}$, we consider the matrix $A=[a_{i,j}]$ the $j$-th
column of which contains the coordinates of the $j$-th base vector
$\xi^{\tau_{1}}_{j}$ of the first cone expressed in terms of the
ordered $L$-basis generating the second cone. The matrix $A$ is an
element of $GL_{n}(\Z)$. The associated monomial mapping $h_{A}$
is defined by $$y_{j}\circ h_{A}:D\rightarrow \R:x\mapsto
\prod_{i=1}^{n}x_{i}^{a_{j,i}}$$ where the domain $D$ consists of
$\R^{n}$ minus the coordinate hyperplanes on which $h_{A}$ is
ill-defined: these are the hyperplanes corresponding to the
variables $x_{k}$ for which not all entries $a_{i,k}$ are
positive.

  It is clear that $h_{A.B}=h_{A}\circ h_{B}$ in points
where both sides are defined. We show that $h_{A^{-1}}$ is defined
on the image of $h_{A}$. Suppose that $x$ belongs to the domain of
$h_{A}$, $y=h_{A}(x)$ and $y_{k}=0$. We have to prove that all
entries in the $k$-th column of $A^{-1}$ are positive, or
equivalently, that $\xi^{\tau_{2}}_{k}$ belongs to $\tau_{1}$. The
fact that $y_{k}=0$ implies the existence of an index $i$ such
that $a_{k,i}>0$ and $a_{j,i}\geq 0$ for all $j$; this means that
$\xi_{i}^{\tau_{1}}$ belongs to $\tau_{2}$. Since $F$ is a fan
$\xi_{i}^{\tau_{1}}$ has to be contained in a common face of
$\tau_{1}$ and $\tau_{2}$, so $a_{j,i}=\delta_{j,k}$. Thus
$\xi_{i}^{\tau_{1}}$ and $\xi^{\tau_{2}}_{k}$ are one and the
same.

  The preceding shows that we have constructed a
well-defined atlas for $X_{L,F}$. Furthermore, the transition
functions map points with positive coordinates in one chart to
points with positive coordinates in another, so the positive part
$X_{L,F}(\R_{+})$ of $X_{L,F}$ is well-defined. When we work with
two fans $F$ and $F'$ we say that $F'$ is finer than $F$
(notation: $F'<F$) if each cone of $F'$ is contained in a cone of
$F$. In this case there exists a natural mapping from $X_{L,F'}$
to $X_{L,F}$: on a standard chart of $X_{L,F'}$ associated to a
$n$-dimensional cone $\tau'$ of $F'$ it is defined as the monomial
mapping associated to the couple of ordered $L$-bases formed by
generators of $\tau'$ and generators of the unique cone $\tau$ of
$F$ containing $\tau'$. Note that the inclusion of $\tau'$ in
$\tau$ implies that the domain of this mapping coincides with the
whole chart. From the nature of this definition it is clear that
all this mappings are compatible as $\tau'$ ranges over the
$n$-dimensional cones of $F'$, so they glue together to a
well-defined analytical mapping $\pi:X_{L,F'}\rightarrow X_{L,F}$.
In the special case where $L=\Z^{n}$ and  $F$ is the positive
orthant $\R_{+}^{n}$ we get a mapping $\pi:X_{L,F'}\rightarrow
\R^{n}$.

  This construction can be generalized by considering two
lattices $L,\,L'$, an $L$-simple fan $F$ and an $L'$-simple fan
$F'$, where $F'<F$. Let $\tau'$ be a $n$-dimensional cone of $F'$
and $\tau$ a cone of $F$ containing $\tau'$. Expressing the
generators of $\tau'$ in the $L$-basis consisting of the ordered
set of generators for $\tau$ yields a monomial map with
nonnegative real exponents, and by gluing we obtain a map $\pi:
X_{L',F'}(\R_{+})\rightarrow X_{L,F}(\R_{+})$.

  The geometric properties of toric varieties are
reflected in the characteristics of the fans used to define them.
The mapping $\pi$ will be proper if and only if the union of the
cones in $F$ coincides with the union of those in $F'$. When
$L=\Z^{n}$ and $F$ is a fan subdividing $\R^{n}_{+}$ and
subordinate to the Newton polyhedron of $f$, the associated
mapping $\pi:X_{\Z^{n},F}\rightarrow \R^{n}$ has a very nice
property: it desingularizes $f$ at the origin of $\R^{n}$. As
always, we assume $f$ to be non-degenerate with respect to its
Newton polyhedron.

  It is clear that the construction of $X_{L,F}$ can be
copied verbatim to the complex case, simply by extending the
transition functions in our atlas to $\C^{n}$, to obtain a complex
analytic variety.

\section{A residue formula}\label{residue}

  Before proceeding, we have to state some conventions.
For every facet $\phi$ of $\Gamma(f)$, let $\xi_{\phi}$ be the
primitive vector (i.e. with components relatively prime in $\N$)
orthogonal to $\phi$. Let $\phi_{1},\ldots \phi_{r}$ be the facets
that contain $\tau_{0}$ and let $\tilde{\tau}_{0}$ be the
$\rho$-dimensional subspace of $\R^{n}$ spanned by these
$\xi_{\phi_{i}}$. Permutating the coordinates of $\R^{n}$ if
necessary, we may assume that $\R^{n}=\tilde{\tau}_{0}\oplus
\sum_{j=\rho+1}^{n} \R e_{j}$ and $\tau_{0}$ is parallel with
$e_{m+1},\ldots, e_{n}$ $(m\geq \rho)$ and with none of the other
$e_{i}$, where $e_{1},\ldots, e_{n}$ is the standard basis for
$\R^{n}$. Putting $N_{i}=l_{\Gamma}(\xi_{\phi_{i}})$ we define $C$
to be the convex hull of
$\{0,\frac{\xi_{\phi_{1}}}{N_{1}},\ldots,\frac{\xi_{\phi_{r}}}{N_{r}},e_{\rho+1},\ldots,e_{n}\}$. If it is not clear from the context which polynomial
$C$ is associated to, we write it explicitly as $C(f)$.

  For $\lambda\in \C$, $\Re(\lambda)>0$ we define a
function $J^{(\lambda)}_{\pm}$ on $D=\{s\in\C\,|\,\Re(s)>0\}$ by
$$J_{\pm}^{(\lambda)}(s)=\int_{\R^{n}_{+}}f_{\pm}(x)^{s}x^{\lambda-1}\varphi(x)\,dx,$$
where $x^{\lambda-1}=x_{1}^{\lambda-1}\cdots x_{n}^{\lambda-1}$.
It is known that $J_{\pm}$, considered as a function in $\lambda$
and $s$, has a meromorphic continuation to the whole of $\C^{2}$
\cite{BerGel}. This essentially comes down to the observation that the
assertion holds when $f$ is a monomial and a reduction to this
particular case via a
resolution of singularities.
If we fix $\lambda$ in $\R_{+}\!\setminus 0$ such that $\lambda
s_{0}\notin \Z$, then the non-integral poles of
$J_{\pm}^{(\lambda)}$ are not greater than $\lambda s_{0}$, and
the polar multiplicity at $\lambda s_{0}$ is at most $\rho$. For
let $F$ be a fan, subdividing the positive orthant, and
subordinate to the Newton polyhedron of $f$, and let
$\pi:X_{\Z^{n},F}\rightarrow \R^{n}$ be the corresponding proper
toric morphism. We know that $\pi$ desingularizes $f$ at the
origin. Using a partition of unity, it suffices to investigate the
poles of the meromorphic continuation of the integral
$$\tilde{J}^{\lambda}_{\pm}(s)=\int_{\R^{n}_{+}}w^{s}\prod_j y_j^{M_js+\kappa_j\lambda-1}\tilde{\varphi}(y)\,dy,$$
which is still defined for $\Re(s)>0$, where $w$ is a nonvanishing
positive analytic function, $\tilde{\varphi}$ is
$\mathcal{C}^{\infty}$ with compact support, and where
$(M_j,\kappa_j)$ is either $(0,1)$, $(1,1)$, or the couple of
numerical data $(N,\nu)$ associated to a ray of $F$. Now, applying
integration by parts, our claim becomes clear.
\begin{theorem}\label{res}
Assume, as always, that $f$ is non-degenerate with respect to its
Newton polyhedron, and furthermore that $\tau_{0}$ is compact; the
latter condition is included only to simplify formulae. When the
support of $\varphi$ is sufficiently small, then
$$\bar{\mu}_{\pm}(\varphi)=n!\,Vol(C)\,\varphi(0)\,PV\,\int_{\R_{+}^{n-\rho}}f_{\tau_{0}}(1,\ldots,1,y_{\rho
+ 1},\ldots,y_{n})^{s_{0}}_{\pm}\,dy\,,$$ where
$$\bar{\mu}_{\pm}=\lim_{s\to
s_0}(s-s_0)^{\rho}\int_{\R_{+}^{n}}f_{\pm}^{s}\varphi\,dx.$$ The
principal value integral $PV\int$ is defined as the value of the
analytic continuation at $\lambda=1$ of the function
$$K_{\pm}(\lambda)= \int_{\R_{+}^{n-\rho}}f_{\tau_{0}}(1,\ldots,1,y_{\rho
+ 1},\ldots,y_{n})^{s_{0}\lambda}_{\pm}y^{\lambda-1}\,dy,$$ where
$K_{\pm}(\lambda)$ is defined for $\lambda\in \R_{+}\!\setminus
0$, and $\lambda s_{0}>-1$. Here $y=\prod_{i=\rho+1}^{n}y_{i}$ and
$dy=dy_{\rho+1}\wedge\ldots\wedge dy_{n}$.
\end{theorem}

\noindent  Some explanation: \\ we will show that
$n!\,Vol(C)\,\varphi(0)\,K_{\pm}(\lambda)$ equals $\lim_{s \to
\lambda s_{0}}(s-\lambda s_{0})^{\rho}J_{\pm}^{(\lambda)}(s)$ on
its domain of definition mentioned above, and that, in particular,
the integral $K_{\pm}$ converges on this domain. This shows that
$K_{\pm}$ has indeed an analytic continuation at $\lambda=1$ -
which is necessarily unique - since we will show, using resolution
of singularities, that, whenever $\Re(\lambda)>0$ and $\lambda
s_{0}\notin \Z$, $(s-\lambda s_{0})^{\rho}J_{\pm}^{(\lambda)}(s)$
is analytic on a neighbourhood of $(s_{0}\lambda,\lambda)$.
Details can be found in the proof.

  In particular, we see that a possible dropping of the
polar multiplicity of $Z_{\pm}(s)$ in $s_{0}$ only depends on
$f_{\tau_{0}}$.
\begin{proof}
%
%
%
%
We may assume that $\dot{\tau}_{0}$ is simple, for the general
case is obtained by subdividing $\dot{\tau}_{0}$ into simple
cones.
Let $L_{1}=\Z^{n}$ and let $F_{1}$ be a $L_{1}$-simple fan,
subordinate to the Newton polyhedron $\Gamma(f)$ of $f$ at 0, and
containing the cone $\dot{\tau_{0}}$ . The natural map
$\pi_{1}:X_{L_{1},F_{1}}\rightarrow \R^{n}$ is an embedded
resolution of singularities of $f$ in a neighbourhood of the
origin in $\R^{n}$ \cite{Arnold}.
  Next, we define the closed submanifold $Y$ of
$X_{L_{1},F_{1}}$, by requiring for every $n$-dimensional cone
$\Delta$ in $F_{1}$ that
$$\begin{array}{l}  U_{L_{1},F_{1},\Delta}\cap Y =  \emptyset\
\mathrm{if} \ \dot{\tau}_{0}\nsubseteq \Delta,$$
 \\[2pt] $$ U_{L_{1},F_{1},\Delta}\cap Y =  \mathrm{locus}\
\{ y_{1}=y_{2}=\ldots=y_{\rho}=0\}\ \mathrm{if}\
\dot{\tau}_{0}\subseteq \Delta,
\end{array}$$
%
%
where $(y_{1},\ldots,y_{n})$ are the standard coordinates in the
chart $U_{L_{1},F_{1},\Delta}$, associated to an ordered basis
$\{\xi_{1},\ldots,\xi_{n}\}$ of $\Delta$ with
$\xi_{1},\ldots,\xi_{\rho}\in \dot{\tau}_{0}$. One can easily
verify that $Y=X_{L_{2},F_{2}}$ where $F_{2}$ is obtained by
projecting the cones in $F_{1}$ containing $\dot{\tau}_{0}$ onto
$\R e_{\rho+1}+\ldots +\R e_{n}=\R^{n-\rho}$, parallel to
$\tilde{\tau}_{0}$, and the lattice $L_{2}$ is the image of
$L_{1}$ under the same projection. Note that the cones of $F_{2}$
are $L_{2}$-simple.
\\[4pt]
Put $L_{3}=\Z e_{\rho+1}+\ldots + \Z e_{n} \subset \R^{n-\rho}$
and let $F_{3}$ be the fan in $\mathbb{R}^{n-\rho}$ induced by all
orthants. Then $X_{L_{3},F_{3}}=(\mathbb{P}_{\R}^{1})^{n-\rho}$,
where $\mathbb{P}_{\mathbb{R}}^{1}$ denotes the real projective
line.

  By refining the fan $F_{1}$ we may suppose that
$F_{2}<F_{3}$. As a consequence of this there exists a natural map
$$\pi_{2}:Y(\mathbb{R}_{+})=X_{L_{2},F_{2}}(\R_{+})\rightarrow
X_{L_{3},F_{3}}(\R_{+})\ .$$ Now the idea is to pull back the
integral defining $K_{\pm}(\lambda)$ along the mapping $\pi_{2}$,
in order to compare $K_{\pm}$ with $\lim_{s\to \lambda
s_{0}}(s-\lambda s_{0})^{\rho}J_{\pm}^{\lambda}(s)$.   Let
$\gamma$ on $(\Pro^{1}_{\R})^{n-\rho}$ be given by
$$\gamma=f_{\tau_{0}}(1,\ldots,1,z_{\rho+1},\ldots,z_{n})_{\pm}^{s_{0}\lambda}\prod_{i=\rho+1}^{n}z_{i}^{\lambda-1}dz_{\rho+1}\wedge\ldots\wedge
dz_{n}\,,$$ where $z_{\rho+1},\ldots, z_{n}$ are standard affine
coordinates on $\R^{n-\rho}$.
With this notation,
$$K_{\pm}(\lambda)=\int_{\R^{n-\rho}_{+}}\gamma=\int_{Y(\R^{+})}\pi_{2}^{*}(\gamma)\,.$$ Let $\Delta$
be a $n$-dimensional cone of $F_{1}$, generated by $\xi_{1},\ldots,
\xi_{n}$, with $\xi_{1},\ldots ,\xi_{\rho}\in \dot{\tau}_{0}$. On
$Y(\R_{+})\cap U_{\Delta}$, where $U_{\Delta}$ is the coordinate
neighbourhood in $X_{L_{1},F_{1}}$ corresponding to $\Delta$, we
have
$$(n!\,Vol(C)\,\varphi(0)\prod_{i=1}^{\rho}N_{i})\pi_{2}^{*}(\gamma)=\frac{(\prod_{i=1}^{\rho}y_{i})\pi_{1}^{*}(\varphi
f_{\pm}^{s_{0}\lambda}x^{\lambda-1}dx)}{dy_{1}\wedge\ldots\wedge
dy_{\rho}}|_{y_{1}=\ldots=y_{\rho}=0}$$ where
$(y_{1},\ldots,y_{n})$ are the standard coordinates associated to
$(\xi_{1},\ldots,\xi_{n})$.
%
%
This equality is straightforward but
crucial for what follows.

  Now we can exploit the special properties of the map
$\pi_{1}$. By \cite{Arnold}, page 202, there exists a neighbourhood of
$\pi_{1}^{-1}(0)$ in $X_{L_{1},F_{1}}$, so that at each point $P$
in this neighbourhood and belonging to $Y(\R_{+})\cap U_{\Delta}$,
we can find a system of local coordinates $y_{1}',\ldots,y_{n}'$
satisfying
\begin{itemize}
\item $y_{i}=y_{i}'$ if $y_{i}(P)=0$\,; in particular this holds
for $i\in \{1,\ldots,\rho\}$, \item
$\pi_{1}^{*}(f^{s}x^{\lambda-1}dx)=v_{1}^{s}v_{2}^{\lambda}v_{3}\prod_{i=1}^{n}y_{i}'^{N_{i}'s+\nu_{i}'-1}dy'$,
where $v_{1},v_{2},v_{3}$ are positive nonvanishing analytic
functions, $(N_{i}',\nu_{i}')=(N_{i},\lambda\nu_{i})$ whenever
$y_{i}(P)=0$, and $(N_{i}',\nu_{i}')$ equals either $(0,1)$ or
$(1,1)$ if $y_{i}(P)\neq 0$, \item the points where $y_{i}\geq 0$
for each $i$ and $f_{\pm}>0$ are exactly the points where
$y'_{j}>0$ for each $j$ satisfying $N'_{j}\neq 0$, and $y'_j\geq
0$ for each $j$ satisfying $y_j(P)=0$.
%
%
%
\end{itemize}
We briefly recall the construction of this new system of local
coordinates $(y'_1,\ldots,y'_n)$. Let us suppose that
$y_1(P)=\ldots=y_s(P)=0$, and that $y_j(P)\neq 0$ when $j>s$. Let
$\gamma$ be the common trace of $\xi_1,\ldots,\xi_s$. By the
definition of $\pi_1$, we can write $f\circ \pi_1$ as
\begin{equation} \label{f0} y_1^{N_1}\ldots
y_s^{N_{s}}(f_0(y_{s+1},\ldots,y_n)+O(y_1,\ldots,y_s))\,.\end{equation}
Now there are two possibilities. If $f_0(y_{s+1},\ldots,y_n)$ is
nonzero at $P$, the factor between brackets is a unit in the local
ring at $P$. If $f_0(y_{s+1},\ldots,y_n)$ vanishes, we can use the
factor between brackets as a new coordinate $y'_{s+1}$, since
$$f_{\gamma}\circ \pi_1=y_1^{N_1}\ldots y_s^{N_{s}}f_0\,,$$
$\pi_1$ induces a local diffeomorphism $\R_0^{n}\rightarrow
\R_0^{n}$, and $f$ is non-degenerate with respect to its Newton
polyhedron. If $P$ were a critical point of $f_0$,
$(1,\ldots,1,y_{s+1}(P),\ldots,y_{n}(P))$ would be a critical
point of $f_{\gamma}$.

It will be important for our purposes, in particular for the
remark following the corollary, that the case
$(N'_i,\nu'_i)=(1,1)$ only occurs when $f_0$ is zero at $P$.

 The choice of local coordinates implies that on $Y(\R_{+})\cap
U_{\Delta}$, the function
$$(n!\,Vol(C)\,\varphi(0)\prod_{i=1}^{\rho}N_{i})\pi_{2}^{*}(\gamma)$$
equals
$$\prod_{i=\rho+1}^{n}y_{i}'^{N_{i}'s_{0}\lambda+\nu_{i}'-1}
(v_{1}^{s_{0}\lambda}v_{2}^{\lambda}v_{3}(\varphi\circ
\pi_{1}))|_{y_{1}'=\ldots=y'_{\rho}=0}\,dy'_{\rho+1}\wedge\ldots\wedge
dy'_{n}\,.$$ Now observe that in the expression
\begin{equation}\label{expr}
\lim_{s\to \lambda s_{0}}(s-\lambda s_{0})^{\rho}
\int_{\R_{+}^{\rho}}\prod_{i=1}^{\rho}y_{i}'^{N_{i}'s+\nu_{i}'-1}
\int_{\R_{+}^{n-\rho}}v_{1}^{s}v_{2}^{\lambda}v_{3}\theta\prod_{i=\rho+1}^{n}y_{i}'^{N_{i}'s+\nu_{i}'-1}dy'\,
\end{equation}
where $\theta$ is a Schwarz function on $\R^{n}$, i.e. a
$\mathcal{C}^{\infty}$-function with compact support, the inner
integral converges for $\lambda>0$ sufficiently small and $s$ near
$s_{0}\lambda$,
since the exponents $N_{i}'s_{0}\lambda+\nu_{i}'-1$, for
$i=\rho+1,\ldots,n$, are either $0$, $\lambda s_{0}$, or
$(N_{i}s_{0}+\nu_{i})\lambda-1>-1$. Hence we can apply the formula
\begin{equation}\label{trick}\lim_{t\searrow\,
t_{0}}(t-t_{0})^{\rho}\int_{[0,a]^{\rho}}\prod_{i=1}^{\rho}z_{i}^{N_{i}(t-t_{0})-1}\psi(t,z)\,dz=\frac{\psi(t_{0},0)}{\prod_{i=1}^{\rho}N_{i}},\end{equation}
which holds for every continuous mapping $\psi$ and any $a\in
\R^{+}_{0}$. Since $N_{i}s_{0}+\nu_{i}=0$ for $i=1,\ldots,\rho$,
this formula yields that \begin{equation}\label{form}\lim_{s\to
\lambda s_{0}}(s-\lambda
s_{0})^{\rho}\int_{X_{L_{1},F_{1}}(\R_{+})}\theta\pi_{1}^{*}(f_{\pm}^{s}x^{\lambda-1}dx)
\end{equation}
is equal to
$$n!\,Vol(C)\int_{\R^{n-\rho}_{+}}\theta(0,\ldots,0,y'_{\rho+1},\ldots,y'_n)\pi_{2}^{*}(\gamma),$$ provided
that the support of $\theta$ is contained in a sufficiently small
neighbourhood of a point of $Y(\R_{+})\cap U_{\Delta}$.

   To conclude the proof of the theorem, one
only has to observe that the expression (\ref{form}) vanishes when
$\theta$ is a Schwarz function with compact support disjoint with
$Y$ (simply apply formula (\ref{trick}) again), and invoke a
suitable partition of unity for $X_{L_{1},F_{1}}$. A similar
construction shows that
$(s-s_{0}\lambda)^{\rho}J^{(\lambda)}_{\pm}(s)$ is analytic in a
neighbourhood of $\{(s_{0}\lambda,\lambda) \in
\C^{2}\,|\,\Re(\lambda)>0,\,\lambda s_{0}\notin \Z\}$\,: simply
apply integration by parts to the integral in (\ref{expr}) with respect to
$y'_{1},\ldots,y'_{\rho}$, and with respect to the $y'_{j}$,
$j>\rho$, with $N'_j\neq 0$, in order to increase their exponent
until its real part becomes greater than $-1$.
\end{proof}

  Note that the compactness of $\tau_{0}$ implies that
for each $i\in\{\rho+1,\ldots,n\}$ there exists an index
$j\in\{1,\ldots,\rho\}$ for which $(\xi_{j})_{i}\neq 0$, so
setting $y_{1},\ldots,y_{\rho}$ equal to zero indeed reduces
$\pi_{1}^{*}(\varphi)$ to $\varphi(0)$. If $\tau_{0}$ fails to be
compact, the factor $\varphi(0)$ has to be replaced by a factor
$\varphi(0,\ldots,0,y_{m+1},\ldots,y_{n})$ in the integrand of the
principal value integral. In particular, the following immediate
consequence of Theorem 1 will still be valid:
\begin{corollary}
  The coefficients $\mu_{\pm}(\varphi)$ for $f$ and
$f_{\tau_{0}}$ differ only by a nonzero factor, which depends only
on the Newton polyhedron of $f$.
\end{corollary}

  \textbf{Remark:} If $f$ satisfies one of the conditions
$\sharp$, then the principal value integral actually converges for
$\lambda=1$, and thus $\mu(\varphi)$ is nonzero, so we recover the
result mentioned in section 3. The fact that condition $(a)$ is
sufficient is obvious. As for condition $(b)$, observe that in
formula (\ref{f0}) in the proof of the theorem, $f_0$ will not
vanish at $P$ if $f\geq 0$ on $\R^{n}$, since this would mean that
$f_0$ has a critical point at $P$. As a consequence, the exponent
$\lambda s_0$ does not occur in (\ref{expr}). To conclude,
condition $(c)$ implies that $f_0$ will not vanish at $P$ in this
case either: since
$$f_{\tau_0}\circ \pi_1=y_1^{N_1}\ldots y_{s}^{N_{s}}(f_0(y_{s+1},\ldots,y_n)+O(y_{1},\ldots,y_s)) \,,$$
the equality $f_0(P)=0$ would induce a zero of $f_{\tau_0}$ in
$\R_0^{n}$.

\section{A new proof of the conjecture}\label{conj}
The objective of this section is to prove the following conjecture
formulated by Denef and Sargos \cite{DeSa2}, Conjecture 3:
\begin{conjecture}
If $\tau_{0}$ is unstable, then $\mu(\varphi)=0$ for any
$\mathcal{C}^{\infty}$-function $\varphi$ with support in a
sufficiently small neighbourhood of 0 in $\R^{n}$.
\end{conjecture}
From the discussion in the last paragraph of section 3, it follows
that we may confine ourselves to the case $s_{0}\notin \Z$.
Moreover, the material in that section shows that is suffices to
prove that $\mu_{\pm}(\varphi)=0$; this is the assertion stated in
Theorem 2. Another - still open - conjecture of Denef and Sargos
claims that the reverse is also true: the vanishing of
$\mu(\varphi)$ whenever the support of $\varphi$ is small enough
implies the instability of $\tau_{0}$.
\begin{theorem}\label{realconj}
If $\tau_{0}$ is unstable with respect to a variable $x_{j}$, then
the polar multiplicity of $Z_{\pm}(s)$ in $s_{0}$ is strictly less
than $\rho$.
\end{theorem}
\begin{proof}
Because of corollary 1, we may suppose that $f=f_{\tau_{0}}$. We
will proceed by constructing an appropriate resolution of
singularities of $f$. To simplify notation we
suppose that $j=n$.

  Let $F$ be a fan subdividing of the positive orthant
$\R_{+}^{n}$ into $\Z^{n}$-simple cones such that $F$ is
subordinate to $\Gamma(f)$. Let $F'$ be the fan consisting of the
simple cones conv($\Delta \cap H_{n}, e_{n}$) with $\Delta \in F$,
where $H_{n}$ is the hyperplane in $\R^{n}$ defined by $x_{n}=0$.
 Let $Y$ be the toric manifold
associated to $(\Z^{n},F')$, and $\pi:Y\rightarrow \R^{n}$ the
natural map.

  We know that $Y$ is nonsingular and $\pi$ is proper.
Furthermore, $\pi$ is an isomorphism on the complement of the
coordinate hyperplanes in $\R^{n}$. In order to prove that $\pi$
is a resolution for $f$ we need to show that locally $f\circ \pi$
and the jacobian $J_{\pi}$ of $\pi$ can be written as the product
of a monomial with a unit.

  Let $\Delta$ be a $n$-dimensional cone of $F'$ spanned
by an ordered basis $\{\xi_{1},\ldots,\xi_{n-1},e_{n}\}$ and $U$
be the associated open part of $Y$. Choosing a point $w$ on $U$ we
may suppose that $1,\ldots,s$ are the indices $i\neq n$ for which
$w_{i}=0$.

  On $U$ the jacobian of $\pi$ is a scalar multiple of
$\prod_{i=1}^{n-1}y_{i}^{\nu_{i}-1}$ and $f\circ\pi$ can be
written as
$$(\prod_{i=1}^{s}y_{i}^{N_{i}})(g(y_{s+1},\ldots, y_{n-1})+h(y_{s+1},\ldots,
y_{n-1})y_{n}+O(y_{1},\ldots,y_{s}))\,.$$ When the $g+h\,y_{n}$
part differs from zero in $w$ we have found our unit; when it
equals zero but $\frac{\partial}{\partial
y_{n}}(g+h\,y_{n})(w)\neq 0$ we can introduce the whole second
factor between brackets as a new variable. So it suffices to show
that $\frac{\partial}{\partial y_{n}}(g+h\,y_{n})=0$ has no
solutions in $\R_{0}^{n-1}\times \R$. But if it has, this means
that $h$ has a zero in $\R_{0}^{n-1}$, contradicting the
definition of unstableness because
$(\prod_{i=1}^{s}y_{i}^{N_{i}})h=f_{\gamma}\circ \pi$, with
$\gamma$ denoting the intersection of the common trace of
$\xi_{1},\ldots,\xi_{s}$ with the hyperplane defined by $x_{n}=1$.

  A resolution of singularities for $f$ determines a set
of candidate poles of $Z_{\pm}(s)$ containing the actual poles,
and provides an upper bound for their polar multiplicities (cf.
\cite{Arnold} and \cite{Var}). This is why we constructed a resolution
that takes into account the instability of $\tau_{0}$: this piece
of extra information yields a sharper upper bound for the polar
multiplicity at $s_{0}$. To make things concrete: under the
assumption that $s_{0}$ is not an integer, this polar multiplicity
is not greater than the maximal number of vectors $\xi_{i}$
occurring as generators of the same cone of the fan $F'$ for which
$N_{i}\neq 0$ and the value $\nu_{i}/N_{i}$ equals $-s_{0}$. As is
easily seen, this condition is equivalent to the property that the
traces of the $\xi_{i}$ contain $\tau_{0}$. Since these $\xi_{i}$
have to be linearly independent and they will automatically be
contained in the hyperplane $x_{n}=0$ (recall that $s_{0}\notin
\Z$), their number can never be greater than $\rho -1$. Now it
becomes clear why we chose this specific form for our resolution:
when we consider a fan subordinate to $\Gamma(f)$, the vectors
$\xi_{i}$ no longer have to be contained in $x_{n}=0$ and their
number can rise up to $\rho$.
\end{proof}

\section{An explicit formula}\label{explicit}
In this section, we give an explicit formula for the residue
$\mu(\varphi)$, in the case where
 $\tau_0$ is a simplex of codimension $1$,
such that the only lattice points in the intersection of $\tau_0$
 with the support of $f$,
 are its vertices. We still suppose that $s_0\notin \Z$, the other case
can be dealt with by introducing a new variable, as was done at the end of
 Section \ref{osc}. Special cases of our explicit
formula were obtained already in the Ph.D. thesis of A. Laeremans
\cite{Laer} under the direction of the first author. We will use
the technique of decoupages,
 developed in \cite{DeSa1}.

Let $f=f_{\tau_0}$ be the polynomial $\sum_{i=1}^{n}\varepsilon_i x^{a_i}$,
 with $x=(x_1,\ldots,x_n)$, with $\varepsilon_i\in \R_0$,
and with $a_i=(a_{i,1},\ldots,a_{i,n})\in \N^n$ for $i=1,\ldots,n$,
 such that the set of vectors $\{a_i\}_i$ linearly independent over $\Q$.
 Here $x^{a_i}$ is short for $\prod_j x_j^{a_{i,j}}$.
We define an $n$-tuple $\gamma$ of positive
real numbers $\gamma_i$, by the expression
$$\vec{1}=(1,\ldots,1)=\sum_{i=1}^{n}\gamma_i a_i\,.$$
Note that, since $t_0\vec{1}$ belongs to $\tau_0$, the $\gamma_i$
satisfy $\sum_i\gamma_i=-s_0$. We denote the real matrix of order
$n$, with the vector $a_i$ as $i$-th column, by $A$. For any
matrix $X$, we will write $X'$ for its transpose. We will denote
the $i$-th column of $A'$ by $a^i$.

Now let $\varphi$ be, as before, a Schwarz function on $\R^{n}$, i.e.
a $\mathcal{C}^{\infty}$-function with compact support $Support(\varphi)$.
We assume that $Support(\varphi)\cap\R^{n}_{+}\subset [0,1]^{n}$.
We denote the linear polynomial $\sum_{i=1}^{n}\varepsilon_iy_i$
by $g(y)$, and we define $\bar{Z}_{\pm}(s)$ by the integral
$$\bar{Z}_{\pm}(s)=\int_{[0,1]^{n}}g_{\pm}(y)^{s}y^{\gamma-1}\varphi(y)dy\,,$$
for $s\in \C$, with $\Re(s)>0$, where $y^{\gamma-1}$ means $\prod_{i=1}^{n} y_i^{\gamma_i-1}$.
 By \cite{DeSa1}, Lemme 3.1,
the function $\bar{Z}_{\pm}(s)$ has a meromorphic continuation to the whole
complex plane, which we will denote again by $\bar{Z}_{\pm}(s)$.

Let $\Sigma$ be a simplicial fan subdividing the positive orthant $\R^{n}_{+}$,
 such that the rays of $\Sigma$
are generated by elements of the set
$\{a^1,\ldots,a^n,e_1,\ldots,e_n\}$, where $(e_i)_{i=1}^{n}$ is
the ordered standard basis for $\R^{n}$. We denote the cones of
$\Sigma$ of dimension $n$ by $\Delta_0,\ldots,\Delta_m$, and we
assume that $\Delta_0$ is generated by the vectors $a^i$. Such a
fan $\Sigma$ always exists, by \cite{DeSa1}, Lemme 2.3.

We will follow the terminology in \cite{DeSa1}.
 We define a function $\mathcal{L}$
by
$$\mathcal{L}:\,]0,1]^{n}\rightarrow \R^{+}_{n}:(x_1,\ldots,x_n)\mapsto
 (-\log x_1,\ldots,-\log x_n)\,.$$
The cones $\Delta_j$ induce a decoupage $\{\mathcal{L}^{-1}(\Delta_j)\}_{j=1}^{m}$ of
$[0,1]^{n}$. It is obvious that $g$ is compatible with this decoupage
 (in the sense of \cite{DeSa1}).
For each $j=1,\ldots,m$, we consider the integral
$$\int_{\mathcal{L}^{-1}(\Delta_j)}g_{\pm}(y)^{s}y^{\gamma-1}\varphi(y)dy,$$
for $s\in \C$, $\Re(s)>0$, and we denote by $Z_{\pm}^{(j)}(s)$ its
meromorphic continuation to the whole of $\C$, which exists by
\cite{DeSa1}, Lemme 5.4.

Now let $\omega$ be the monomial transformation associated to $\Delta_0$, i.e.
$$\omega :[0,1]^{n}\rightarrow [0,1]^{n} :(x_i)_i\mapsto (\prod_{i=1}^{n}x_i^{a_{i,k}})_k\,. $$
If we pull back the integral defining $Z^{(0)}_{\pm}$ along $\omega$, we get
$$Z^{(0)}_{\pm}(s)=|det\,A|\int_{[0,1]^{n}}f_{\pm}(x)^{s}\varphi\circ \omega(x) dx\,.$$

\begin{lemma}
The function $Z^{(0)}_{\pm}(s)$ is the only term in the sum
$\bar{Z}_{\pm}(s)=\sum_{j=0}^{m}Z^{(j)}_{\pm}(s)$ with a pole at $s=s_0$.
\end{lemma}

\textit{Remark.} Our explicit computations below, will show that $\bar{Z}_{\pm}(s)$ has,
 indeed, a pole at $s=s_0$.

\begin{proof}
First, we show that the vector $\vec{1}$ is contained in the interior of $\Delta_0$.
We define an $n$-tuple of real numbers $\zeta=(\zeta_i)_i$ by $\zeta A'=\vec{1}$.
This means that $A\zeta'=\vec{1}'$, hence $\sum \zeta_iz_i=0$ is the equation of $\tau_0$,
which implies that $\zeta_i> 0$ for every $i$.

We subdivide $\Delta_0$ into $n$ simplicial subcones
$\Delta_{0,l}$, where $\Delta_{0,l}$ is generated by $\vec{1}$,
and $\{a^1,\ldots,\hat{a^l},\ldots,a^n\}$. Let $\Sigma'$ be the
induced subdivision of the fan $\Sigma$. We say that a ray of
$\Sigma'$, generated by some vector $\xi$,
 contributes to the pole $s_0$ of $\bar{Z}{\pm}$, if there exists
a positive integer $\alpha$, such that
$$s_0= -\frac{\!<\!\xi,\gamma\!>\!+\alpha}{l_{\Gamma(g)}(\xi)}\,,$$
where $l_{\Gamma(g)}$ is, as in Section \ref{poly},
 the trace map associated to the Newton polyhedron $\Gamma(g)$ of $g$.
If $l_{\Gamma(g)}(\xi)$ is zero, we define the right part of the equation to be
 equal to $-\infty$.
It is clear that the vector $\vec{1}$ contributes to $s_0$,
 and we will show it is the only ray of $\Sigma'$ that does.
By \cite{DeSa1}, Lemme 5.4, this concludes the proof.

For a tuple $\xi=(\xi_i)$ of real numbers, we will write $\min(\xi)$ for $\min_i\{\xi_i\}$.
Let $T$ be the torus $\R^{n}_{0}$.
For $\xi$ in $T(\R_{+})$, and $\xi\notin \R\vec{1}$, we see that
\begin{eqnarray*}
\frac{\!<\!\xi,\gamma\!>\!}{l_{\Gamma(g)}(\xi)}
&=&\frac{\!<\!\xi,\gamma\!>\!}{\min(\xi)}
\\ &>&\frac{\!<\!\min(\xi)\vec{1},\gamma\!>\!}{\min(\xi)}
\\ &=& -s_0
\end{eqnarray*}
This shows that $\vec{1}$ is the only vector contibuting to $s_0$.
\end{proof}

Hence, in order to know the residue of $Z^{(0)}_{\pm}(s)$ at $s=s_0$,
 it suffices to
study the residue of  $\bar{Z}_{\pm}(s)$ at $s_0$.

The function $Z^{(0)}_{\pm}(s)$ is more or less
 the meromorphic function whose residue in $s_0$
we want to investigate. The problem is that the
 function $(\varphi\circ\omega)\chi_{[0,1]^{n}}$
is not $\mathcal{C}^{\infty}$, and that the support of $\varphi\circ\omega$
cannot be chosen in an arbitrarily small neighbourhood of $0$,
 since it contains $\omega^{-1}(0)$ as soon as $\varphi(0)\neq 0$.
 The following lemma deals with these difficulties.

\begin{lemma}
Let $\{\psi_{\alpha}\}$ be a partition of unity for $\R^{n}$,
with $\psi_{\alpha_0}\equiv 1$ on a neighbourhood of $0$.
 If the supports of $\varphi$ and $\psi_{\alpha_0}$ are sufficiently small,
$$\lim_{s\to s_0}(s-s_0)^{\rho}Z^{(0)}_{\pm}(s)=
|det\,A|\lim_{s\to s_0}(s-s_0)^{\rho}Z_{\alpha_0}(s),$$
where we write $Z_{\alpha_0}$ for
 the meromorphic continuation of
$$\int_{\R_+^{n}}f_{\pm}(x)^{s}
\varphi\circ \omega(x)\psi_{\alpha_0}(x) dx .$$
\end{lemma}
\begin{proof}
Our proof is similar to the proof of the previous lemma.
We can construct a subcone $\Delta'$ of $\Delta_0$, containing
$\vec{1}$ in its interior, such that
$$\omega^{-1}(\mathcal{L}^{-1}(\Delta')\cap
 Support(\varphi))
\subset \psi_{\alpha_0}^{-1}(1),$$
provided the support of $\varphi$ is sufficiently small.
If we extend $\{\Delta'\}$ to a subdivision of $\Delta_0$,
 the only cone contributing to the residue of $\bar{Z}_{\pm}(s)$ at $s=s_0$
 will be $\Delta'$ itself.

We denote by $Z_{\pm}^{\Delta'}$ the meromorphic continuation of
$$\int_{(\mathcal{L}\circ \omega)^{-1}(\Delta')}f_{\pm}(x)^{s}
\varphi\circ \omega(x)dx .$$
We introduce a new function $J_{\Delta'}^{(\lambda)}(s)$, which is defined,
for $\Re (\lambda)>0$ and $\Re (s)>0$, by
$$J_{\Delta'}^{(\lambda)}(s)=\int_{\mathcal{L}^{-1}(\Delta')}
g_{\pm}(y)^{s\lambda}
y^{\lambda\gamma-1}
\varphi(y)dy \,.$$
Pulling back the integral via $\omega$ yields
$$\int_{(\mathcal{L}\circ \omega)^{-1}(\Delta')}
f_{\pm}(x)^{s\lambda}x^{\lambda-1}
\varphi\circ \omega(x) dx .$$

Slightly adapting the proofs of \cite{DeSa1}, Lemme 3.1 and Lemme 5.4,
we see that this function has a meromorphic continuation to $\C^{2}$,
which we will again denote by $J_{\Delta'}^{(\lambda)}(s)$.
If we fix  $\lambda$ in $\R_{+}\!\setminus 0$ such that $\lambda
s_{0}\notin \Z$, then the non-integral poles of
$J_{\Delta'}^{(\lambda)}$ are not greater than $\lambda s_{0}$, and
the polar multiplicity at $\lambda s_{0}$ is at most $\rho$.
Similarly, we define a function $J_{\alpha_0}^{(\lambda)}(s)$ as
the meromorphic continuation of
$$\int_{\R_+^{n}}f_{\pm}(x)^{s\lambda}x^{\lambda-1}
\varphi\circ \omega(x)\psi_{\alpha_0}(x) dx .$$

 Now observe that formula \ref{trick}
 in our proof of Theorem \ref{res} holds, as soon as $\psi$ is continuous
in a neighbourhood of $(t_0,0)$. This implies that, when the support of
$\psi_{\alpha_0}$ is sufficiently small, and $\lambda>0$ is sufficiently small,
\begin{equation}\label{eq}
\lim_{s\to \lambda s_0}(s-\lambda s_0)^{\rho}J_{\Delta'}^{(\lambda)}(s)=
|det\,A|\lim_{s\to \lambda s_0}(s-\lambda s_0)^{\rho}J_{\alpha_0}^{(\lambda)}(s).
\end{equation}

The function $(s-\lambda s_0)^{\rho}J_{\Delta'}^{(\lambda)}(s)$ is meromorphic
on $\C^{2}$, and by (\ref{eq}), the plane $s=\lambda s_0$ is not
contained in its polar locus (this also
follows from the results in \cite{DeSa1}).
 Hence, its restriction to this plane is
a meromorphic function in $\lambda$, and by the identity principle,
it coincides with the restriction of
$|det\,A|(s-\lambda s_0)^{\rho}J_{\alpha_0}^{(\lambda)}(s).$
Taking values in $\lambda=1$ yields
\begin{equation*}
\lim_{s\to s_0}(s-s_0)^{\rho}Z_{\pm}^{\Delta'}(s) =
|det\,A|\lim_{s\to s_0}(s-s_0)^{\rho}Z_{\alpha_0}(s).
\end{equation*}

Hence,
$$\lim_{s\to s_0}(s-s_0)^{\rho}Z^{(0)}_{\pm}(s)=
|det\,A|\lim_{s\to s_0}(s-s_0)^{\rho}Z_{\alpha_0}(s),$$
since both sides are equal to
the residue of  $\bar{Z}_{\pm}(s)$ at $s_0$.
\end{proof}

 It follows from Theorem \ref{res}, that
$$\lim_{s\to s_0}(s-s_0)^{\rho}Z^{(0)}_{\pm}(s)=
\lim_{s\to s_0}(s-s_0)^{\rho}\tilde{Z}(s),$$
with $\tilde{Z}$ the meromorphic continuation of
$$\int_{[0,1]^{n}}f_{\pm}(x)^{s}
\phi(x) dx, $$
 whenever $\phi$ is a Schwarz function with
sufficiently small support, and $\phi(0)=\varphi(0)$.

We define $\bar{I}(t)$ as the oscillating integral
$$\bar{I}(t) =\int_{[0,1]^{n}}e^{itg(y)}y^{\gamma-1}\varphi(y)dy\,. $$
As the coefficient of $t^{s_0}$
 in its asymptotic expansion only depends on $\varphi(0)$,
provided the support of $\varphi$ is sufficiently small, we may assume that
$\varphi$ is of the form $\prod_{i=1}^{n}\theta_i(y_i)$, where $\theta_i$
is a Schwarz function on $\R$.
Then $\bar{I}(t)$ splits into a product of integrals of the type
$$\mathcal{I}(t;\varepsilon,\eta) = \int_{0}^{1}e^{it\varepsilon z}z^{\eta-1}\theta(z)dz\,,$$
where $\varepsilon$ and $\eta$ are non-zero real numbers, $\eta>0$.
By \cite{Arnold}, 7.2.3 (11), we get, as $t\to +\infty$,
$$\mathcal{I}(t;\varepsilon,\eta) \sim \theta(0) \frac{\Gamma(\eta)}{(-i\varepsilon t)^{\eta}},$$
where $\arg(\pm it)=\pm\pi/2$, and $\Gamma$ is the Gamma function.
Bringing these factors together, we see that
the coefficient of $t^{s_0}$ in the asymptotic expansion
of $\bar{I}(t)$ is equal to
$$\varphi(0)(\prod_{j=1}^{n} \Gamma(\gamma_j)|\varepsilon_j|^{-\gamma_j}
e^{sign(\varepsilon_j)\frac{i\pi}{2}\gamma_j}). $$

\begin{theorem}\label{formula}
Suppose that  $f$ is non-degenerate
with respect to its Newton polyhedron,
 and $\tau_0$ is a simplex of codimension $\rho=1$,
such that the only lattice points in the intersection of $\tau_0$
 with the support of $f$,
 are its vertices. We suppose as well that $s_0\notin\Z$. If $\varphi$
is a Schwarz function with sufficiently small support,
$$\mu(\varphi)=
|det\,A|^{-1}\varphi(0)(\prod_{j=1}^{n}\Gamma(\gamma_j)
|\varepsilon_j|^{-\gamma_j})\sum_{\beta\in\{-1,1\}^{n}}
\prod_{j=1}^{n}e^{sign(\varepsilon_j\beta^{a_j})\frac{\pi i}{2}\gamma_j} .$$
\end{theorem}
\begin{proof}
This follows immediately from our residue formula, and the
arguments above.
\end{proof}


Conjecture 4 of \cite{DeSa2} predicts the following remarkable combinatorial assertion.
\begin{conjecture}
With the notation and hypotheses of Theorem \ref{formula},
if $\tau_0$ is stable, then
$$\sum_{\beta\in\{-1,1\}^{n}}
\prod_{j=1}^{n}e^{sign(\varepsilon_j\beta^{a_j})\frac{\pi i}{2}\gamma_j}\neq 0.$$
\end{conjecture}
\section{A complex residue formula}\label{sectioncomplex}
In this section, we will establish te complex analogue of our residue
formula in Theorem \ref{res}.
 The complex local zeta function $Z(s)$, associated
to a complex analytic germ $f$ at $(\C^{n},0)$, is defined for
$s\in \C,\,\Re(s)>0$ as
$$\int_{\C^{n}}|f|^{2s}\varphi(z)dz$$ where $\varphi$ is a
positive $\mathcal{C}^{\infty}$-function with sufficiently small
support, and the integration is conducted with respect to the
Haar-measure on $\C^{n}$ (this is just the Lebesgue measure on
$\C^{n}=\R^{2n}$). In the definition, we use $|f|^{2s}$,
rather than $|f|^{s}$, for reasons of uniformity: $|.|^{2}$ is the
modulus associated to the Haar-measure on the local field $\C$, as
are $|.|$ in the real, and $|.|_{p}$ in the $p$-adic case \cite{Igusa:intro}.
 It is
known that $Z$ has a meromorphic continuation to the whole of
$\C$, which we will denote again by $Z$. We again assume that
$f$ vanishes at the origin, and that
$f$ is non-degenerate with respect to its Newton polyhedron.
We define $\tau_0$, $s_0$, and $\rho$, as in the real case.
We will use the same notation and conventions as in Section \ref{residue}.

As before, a toric embedded resolution establishes $s_0$ as the
largest non-trivial candidate pole of $Z(s)$, where we say that a
pole is trivial if it is integer and has multiplicity $1$. We will
assume that $s_0\not\in \Z$. In this case, the polar multiplicity
of $s_0$ is at most $\rho$. We want to determine the residue
$$|\mu|(\varphi)=\lim_{s\to s_0}(s-s_0)^{\rho}Z(s)\,.$$

\noindent For $\lambda\in \C$, $\Re(\lambda)>0$ we define a
function $J^{(\lambda)}$ on $D=\{s\in\C\,|\,\Re(s)>0\}$ by
$$J^{(\lambda)}(s)=\int_{\C^{n}}|f(x)|^{2s}|x|^{2\lambda-2}\varphi(x)\,dx,$$
where $|x|^{2\lambda-2}=|x_{1}|^{2\lambda-2}\cdots
|x_{n}|^{2\lambda-2}$. As in the real
case, $J$, considered as a
function in $\lambda$ and $s$, has a meromorphic continuation to
the whole of $\C^{2}$, and
 if we
fix $\lambda$ in $\R_{+}\!\setminus 0$ such that $\lambda
s_{0}\notin \Z$, then the non-integral poles of $J^{(\lambda)}$
are not greater than $\lambda s_{0}$, and the polar multiplicity
at $\lambda s_{0}$ is at most $\rho$.

The following theorem is the complex counterpart of Theorem \ref{res}.
\begin{theorem}\label{rescomp}
Assume, as always, that $f$ is non-degenerate with respect to its
Newton polyhedron, and furthermore that $\tau_{0}$ is compact; the
latter condition is included only to simplify formulae. When the
support of $\varphi$ is sufficiently small, then
$$|\mu|({\varphi})=\pi^{\rho}n!\,Vol(C)\,\varphi(0)\,PV\,\int_{\C^{n-\rho}}|f_{\tau_{0}}(1,\ldots,1,z_{\rho
+ 1},\ldots,z_{n})|^{2s_{0}}\,dz\,,$$ where the principal value
integral is defined as the value of the meromorphic continuation
at $\lambda=1$ of the function
$$K(\lambda)= \int_{\C^{n-\rho}}|f_{\tau_{0}}(1,\ldots,1,z_{\rho
+ 1},\ldots,z_{n})|^{2s_{0}\lambda}|z|^{2\lambda-2}\,dz,$$ where
$K(\lambda)$ is defined for $\lambda\in \R_{+}\!\setminus 0$ and
$\lambda s_{0}>-1$. Here $z=\prod_{i=\rho+1}^{n}z_{i}$, and $dz$ is
the Haar- measure on $\C^{n-\rho}$.
\end{theorem}

\noindent Some explanation: \\ we will show, again, that
$\pi^{\rho}n!\,Vol(C)\,\varphi(0)\,K(\lambda)$ equals $\lim_{s \to
\lambda s_{0}}(s-\lambda s_{0})^{\rho}J^{(\lambda)}(s)$ on its
domain of definition mentioned above,
 and that, in particular,
the integral $K$ converges on this domain.
This shows that
$K$ has indeed an analytic continuation at $\lambda=1$ -
which is necessarily unique - since we will show, using resolution
of singularities, that, whenever $\Re(\lambda)>0$ and $\lambda
s_{0}\notin \Z$, $(s-\lambda s_{0})^{\rho}J^{(\lambda)}(s)$
is analytic on a neighbourhood of $(s_{0}\lambda,\lambda)$.
 Details can be found in the
proof.

\noindent In particular, we see that the dropping of the polar
multiplicity of $Z(s)$ in $s_{0}$ only depends on $f_{\tau_{0}}$.

\begin{proof}
The proof is almost the same as in the real case.
We define lattices $L_i$ and fans $F_i$, $i=1,\ldots,3$, as before, as well
as the resolution morphism $\pi_1$, and the submanifold $Y$ of $X_{L_1,F_1}$.
Refining the fan $F_{1}$, we may suppose that
$F_{2}<F_{3}$. We define $L_3'$ as $\frac{1}{N}L_3$, with $N\in \N_0$.
The toric morphism $$\pi':X_{L'_3,F_3}\rightarrow X_{L_3,F_3}$$ is
an $N(n-\rho)$-fold cover of $(\Pro^{1}_{\C})^{n-\rho}$,
 ramified over the orbits of codimension one. For an appropriate choice of $N$,
a set of base vectors of $L_2$ has integer coordinates with respect to $L'_3$.
 As a consequence of this, there exist natural maps
$$\pi_{2}:X_{L'_{3},F_{2}}\rightarrow
X_{L_{3},F_{3}}\ \mathrm{and}\ \pi_{3}:X_{L'_{3},F_{2}}\rightarrow
X_{L_{2},F_{2}}=Y\ .$$

Let $\gamma$ on
$(\Pro^{1}_{\C})^{n-\rho}$ be given by
$$\gamma=|f_{\tau_{0}}(1,\ldots,1,w_{\rho+1},\ldots,w_{n})|^{2 s_{0}\lambda}\prod_{i=\rho+1}^{n}|w_{i}|^{2\lambda-2}dw_{\rho+1}\wedge\ldots\wedge
dw_{n}\,,$$ where $w_{\rho},\ldots, w_{n}$ are standard affine
coordinates on $\C^{n-\rho}$.
 With this notation,
$$K(\lambda)=\int_{\C^{n-\rho}}\gamma=\int_{X_{L_{3},F_{2}}}\pi_{2}^{*}(\gamma)\,.$$ Let $\Delta$
be a $n$-dimensional cone of $F_{1}$, generated by $\xi_{1},\ldots,
\xi_{n}$, with $\xi_{1},\ldots ,\xi_{\rho}\in \dot{\tau}_{0}$. On
$Y\cap U_{\Delta}$, where $U_{\Delta}$ is the coordinate
neighbourhood in $X_{L_{1},F_{1}}$ corresponding to $\Delta$, we
have
$$(n!\,Vol(C)\,\varphi(0)\prod_{i=1}^{\rho}N_{i})\pi_{2}^{*}(\gamma)=\pi_{3}^{*}(\frac{(\prod_{i=1}^{\rho}|y_{i}|^{2})\pi_{1}^{*}(\varphi
|f|^{2
s_{0}\lambda}|z|^{2\lambda-2}dz)}{dy_{1}\wedge\ldots\wedge
dy_{\rho}}|_{y_{1}=\ldots=y_{\rho}=0})$$ where
$(y_{1},\ldots,y_{n})$ are the standard coordinates associated to
$(\xi_{1},\ldots,\xi_{n})$.

We can again find local coordinates $y_i'$, around each point $P$ of
$Y\cap U_{\Delta}$ in a neighbourhood of $\pi_1^{-1}(0)$, such that
\begin{itemize}
\item $y_{i}=y_{i}'$ if $y_{i}(P)=0$\,; in particular this holds
for $i\in \{1,\ldots,\rho\}$, \item
$\pi_{1}^{*}(|f|^{2s}|z|^{2\lambda-2}dz)=|v_{1}|^{s}|v_{2}|^{\lambda}|v_{3}|
\prod_{i=1}^{n}|y_{i}'|^{2N_{i}'s+2\nu_{i}'-2}dy'$,
where $v_{1},v_{2},v_{3}$ are nonvanishing analytic
functions, $(N_{i}',\nu_{i}')=(N_{i},\lambda\nu_{i})$ whenever
$y_{i}(P)=0$, and $(N_{i}',\nu_{i}')$ equals either $(0,1)$ or
$(1,1)$ if $y_{i}(P)\neq 0$.
\end{itemize}

In the expression
\begin{equation}
(s-\lambda s_{0})^{\rho}
\int_{\C^{\rho}}\prod_{i=1}^{\rho}|y_{i}'|^{2 N_{i}'s+2\nu_{i}'-2}
\int_{\C^{n-\rho}}|v_{1}|^{s}|v_{2}|^{\lambda}|v_{3}|\theta\prod_{i=\rho+1}^{n}|y_{i}'|^{2
N_{i}'s+2\nu_{i}'-2}dy'\,
\end{equation}
where $\theta$ is a Schwarz function on $\C^{n}=\R^{2n}$, the
inner integral converges for $\lambda$ sufficiently small and
$s>s_{0}\lambda$,
%
%
since the exponents $N_{i}'s_{0}\lambda+\nu_{i}'-1$, for
$i=\rho+1,\ldots,n$, are either $0$, $\lambda s_{0}$, or
$(N_{i}s_{0}+\nu_{i})\lambda-1>-1$. Hence we can apply the formula
$$\lim_{t\searrow\,
t_{0}}(t-t_{0})^{\rho}\int_{\|z\|\leq 1
}\prod_{i=1}^{\rho}z_{i}^{2
N_{i}(t-t_{0})-2}\psi(t,z)\,dz=\frac{\pi^{\rho}\psi(t_{0},0)}{\prod_{i=1}^{\rho}N_{i}},$$
which holds for every continuous mapping $\psi$ on $\R\times
\C^{\rho}$.

\noindent To conclude the proof of the theorem, one only has to
observe that $$\lim_{s\to \lambda s_{0}}(s-\lambda
s_{0})^{\rho}\int_{X_{L_{1},F_{1}}}\theta\pi_{1}^{*}(|f|^{2s}|z|^{2\lambda-2}dz)=0$$
when $\theta$ is a Schwarz function with compact support disjoint
with $Y$, and invoke a suitable partition of unity for
$X_{L_{1},F_{1}}$. A similar construction shows that
the function $(s-s_{0}\lambda)^{\rho}J^{(\lambda)}(s)$ is analytic in a
neighbourhood of $$\{(s_{0}\lambda,\lambda) \in
\C^{2}\,|\,\Re(\lambda)>0,\,\lambda s_{0}\notin \Z\}.$$
\end{proof}

\noindent
 If $\tau_{0}$ fails to be compact, the factor $\varphi(0)$ has
to be replaced by a factor
\\$\varphi(0,\ldots,0,z_{m+1},\ldots,z_{n})$ in the integrand of
the principal value integral ($m$ is defined as in the real case).
 In particular, the following
immediate consequence of Theorem \ref{rescomp} will still be valid:
\begin{corollary}
The coefficients $|\mu|(\varphi)$ for $f$ and $f_{\tau_{0}}$
differ only by a nonzero factor, which depends only on the Newton
polyhedron of $f$ .
\end{corollary}
If $s_{0}>-1$, it follows from the proof that the principal value
integral actually converges for $\lambda=1$, and hence that
$|\mu|(\varphi)$ is nonzero.
\section{The stability conjecture}\label{sectionstability}
Denef and Sargos stated their conjecture in the complex case, as well
 \cite{DeSa2}. As always, we suppose that $s_0\notin \Z$.
\begin{conjecture}[Denef-Sargos]
Let $f$ be a singular complex analytic germ at $0\in \C^{n}$, with $f(0)=0$,
and let $\varphi$ be a Schwarz function on $\C^{n}$
 with sufficiently small support, such that $\varphi(0)\neq 0$.
We denote by $Z(s)$ the complex local zeta function associated to $f$ and
$\varphi$, and we suppose $s_0\notin\Z$.
The polar multiplicity of $Z(s)$ at $s_{0}$ equals
$\rho$ if and only if $\tau_{0}$ is a stable face.
\end{conjecture}

We can copy the proof of Theorem \ref{realconj} verbatim to obtain a proof
of its complex analogue.

\begin{theorem}\label{cconj}
If $\tau_{0}$ is unstable with respect to a variable $z_{j}$, then
the polar multiplicity of $Z(s)$ in $s_{0}$ is strictly less
than $\rho$.
\end{theorem}

However, in the complex case, we can say more.
We start by proving the following statement.

\begin{prop}\label{complexconj}
If $\tau_{0}$ is compact, and
$\rho=1$, the complex local zeta function $Z(s)$
has a pole of order $\rho$ at $s=s_0$.
\end{prop}

\begin{proof}
 First, we reduce to the case where $f$ has an isolated singularity at
the origin.
If not, we can
always modify $f$ without changing $f_{\tau_{0}(f)}$, which is the
only part of $f$ that really matters as far as the polar
properties of $s_{0}$ are concerned, by adding, for each $i=1,\ldots,n$,
 a monomial
$b_{i}x_{i}^{a_{i}}$ to $f$, with $a_{i}$
sufficiently large. 
The new polynomial $f$, obtained in this way, will still be
non-degenerate with respect to its Newton polyhedron, and
furthermore, it will have an isolated singularity at the origin,
provided we make suitable choices for $a_{i}, b_{i}$ (cf.
\cite{Kou}).

Loeser showed in \cite{Loes}, that the fact that $f$ has an isolated
singularity at the origin, implies that
the set $\{-(\alpha+m)\,|\,\alpha\in Sp(f),\,m\in\N_0\}$
is contained in the set of poles of $Z(s)$, where
$Sp(f)$ is the spectrum of Steenbrink-Varchenko associated to $f$.
Since $-s_{0}-1$ is the smallest value in $Sp(f)$ (cf. \cite{Kuli})
we can conclude that $s_{0}$ is indeed a pole of $Z$.
\end{proof}

This Proposition does not contradict the Conjecture,
 since in the case $\rho=1$,
$\tau_0$ will automatically be stable. Using this Proposition, and our
complex residue formula, we can derive a combinatorial sufficient condition
that implies the Conjecture.

\begin{lemma}\label{lift}
Suppose that $\tau_0(f)$ is compact, and
 that we can find a complex analytic function $g(z_{\rho},\ldots,z_n)$
on a neighbourhood of the origin in $\C^{n-\rho+1}$, such that
\begin{itemize}
\item $g(0)=0$ and $g$ has a singularity at the origin
\item $g$ is non-degenerate with respect to its Newton polyhedron
\item  $\tau_{0}(g)$ is compact, and has codimension 1
\item$s_0(g)=s_0(f)$
\item $\tilde{\tau}_{0}(g)\oplus \sum_{j=\rho+1}^{n}\R e_{j}=\R^{n-\rho+1}$.
\item $g_{\tau_0(g)}(1,z_{\rho+1},\ldots,z_n)=
f_{\tau_0(f)}(1,\ldots,1,z_{\rho+1},\ldots,z_n)$
\end{itemize}
Then $Z(s)$ has a pole at $s=s_0$ of order $\rho$.
\end{lemma}
\begin{proof}
We will denote the complex local zeta function, associated to a complex
analytic germ $h$, and a Schwarz function $\phi$, by $Z_{h,\phi}(s)$.
 We deduce from Proposition \ref{complexconj} and our residue formula that
\begin{eqnarray*}
\lim_{s\to
s_{0}}(s-s_{0})^{\rho}Z_{f,\varphi}(s)&=
&c\lim_{s\to s_{0}}(s-s_{0})Z_{g,\tilde{\varphi}}(s)
\\&\neq&0
\end{eqnarray*}

\noindent where
$\tilde{\varphi}(x_{\rho},\ldots,x_{n})=\varphi(0,\ldots,0,x_{\rho},\ldots,x_{n})$
and $c$ is a nonzero constant. So $s_{0}$ is a pole of order
$\rho$ of $Z(s)$.
\end{proof}

The natural choice for the function $g$ would be
$$g(z_{\rho},\ldots,z_n)=f_{\tau_{0}}(1,\ldots,1,z_{\rho},\ldots,z_n),$$
maybe after permutating the first $\rho$ variables $z_i$.
However, this does not work in general, because it can happen that
$g\neq g_{\tau_0(g)}$, as is illustrated by the following example.

\begin{example}
Consider the polynomial function $$f(z_1,z_2,z_3,z_4) =
  z_2^{2}z_3^{3}z_4 + z_1^{2}z_3z_4^{3} + z_1^{2}z_2^{2}z_3z_4.$$
It is easy to see that $\tau_0(f)=supp(f)$, that $\rho=2$, and
that $\tau_0$ is stable. The faces of $\Gamma(f)$ that contain $\tau_0$
are contained in the hyperplanes $z_2+z_4=3$, and $z_1+z_3=3$,
 which means that any coordinate could figure as $z_1$. However,
no matter how we permutate the coordinates,
$g(z_2,z_3,z_4)=f(1,z_2,z_3,z_4)$ will never satisfy $g_{\tau_0(g)}=g$.
\end{example}

However, this choice does work under
some additional assumptions. We denote by $\pi_{j}$ the projection
of $\R^n$ onto $\sum_{k=j}^{n}\R e_j$, with $1\leq j\leq n$.

\begin{theorem}\label{assum}
Suppose that $\tau_0$ is stable and compact,
and suppose that, maybe after permutating $(z_1,\ldots,z_{\rho})$,
the function $g$ defined by
$$g(z_{\rho},\ldots,z_n)=f_{\tau_{0}}(1,\ldots,1,z_{\rho},\ldots,z_n)$$
satisfies $\tau_0(g)=\pi_{\rho}(\tau_0(f))$.
 Then $Z(s)$ has a pole of order $\rho$
at $s=s_0$.
\end{theorem}
\begin{proof}
It suffices to check the conditions of Lemma \ref{lift}.
By assumption, $g_{\tau_0(g)}=g$, and $\tau_0(g)$ is compact.

Let us first show that
the dimension of $\tau_0(g)$ is still equal to $n-\rho$,
 hence $\tau_0(g)$ has codimension $1$.
Our assumption implies that
$$\pi_{m}(\tau_0(f))=\tau_0(f(1,\ldots,1,z_m,\ldots,z_n)),$$
for each $m=1,\ldots,\rho$.
Choose vertices $v_{0},\ldots,v_{n-\rho}$ of $\tau_0(f)$, spanning
the affine subspace $V$ of $\R^{n}$ generated by $\tau_0(f)$.
 Since $\tau_0(f)$ is compact,
$dim\,\pi_2(V)=dim\,V$, because $V$ can not contain vectors parallel
to one of the coordinate axes. Repeating this argument, shows that
 $\tau_0(g)$ has codimension $1$. It is clear that
$\tilde{\tau}_0(g)\oplus\sum_{j=\rho+1}^{n}\R e_j=\R^{n-\rho+1}$.

 The stability of $\tau_{0}(g)$
follows from the stability of $\tau_{0}(f)$. For suppose that
$\tau_{0}(g)$ is unstable with respect to $z_{k}$. This means that
$\tau_{0}$ is a pyramid with top in $z_{k}=1$ and base in
$z_{k}=0$ \cite{DeSa2}. Then the same holds for the face
$\tau_{0}(f)$, since it cannot
 be
parallel to one of the projection axes. This is impossible, as we assumed that
$\tau_0(f)$ is stable.
The stability of $\tau_0(g)$ implies, in particular, that the origin is a singularity
 for $g$.

 To conclude, we show that $g$ is non-degenerate with respect to
its Newton polyhedron. Let $\sigma$ be a face of $\tau_{0}(g)$,
and suppose that $(\alpha_{\rho},\ldots, \alpha_{n})$ is a
singular point of $g_{\sigma}$. The face $\sigma$ corresponds to a
face $\sigma'$ of $\tau_{0}(f)$. Our aim is to prove that one of
the $\alpha_{i}$ is zero, and we proceed by showing that
$\alpha=(1,\ldots,1,\alpha_{\rho},\ldots, \alpha_{n})$ is a
singular point of $f_{\sigma'}$. The quasihomogeneity of
$f_{\sigma'}$ implies that
$N_jf_{\sigma'}=\sum_{i=1}^{n}(\xi_{\phi_j})_{i}z_{i}\frac{\partial
f_{\sigma'}}{\partial z_{i}}$ for $j=1,\ldots,r$.
Here, the $\phi_j$ are the facets of $\Gamma(f)$ containing $\tau$,
$\xi_{\phi_j}$ is the primitive normal vector on $\phi_j$, and $N_j$ equals
$l_{\Gamma(f)}(\xi_{\phi_j})$.
 Applying both
sides of the expression to $\alpha$,
 yields a linear system of equations of rank
$\rho-1$ in the first $\rho-1$ partial derivatives of $f$ in
$\alpha$, obliging all of these to become zero as well.
\end{proof}

In spite of the additional assumptions, the proof of Theorem \ref{assum}
gives some intuition about the reasons behind the stability condition:
it guarantees that $g$ still has a singularity at the origin.

%

\begin{corollary}\label{cor3}
If $\tau_0$ is a stable, compact face of dimension $1$, then
$Z(s)$ has a pole of order $\rho$ at $s=s_0$.
\end{corollary}
\begin{proof}
By Proposition \ref{complexconj}, we may assume $n\geq 3$. Let
$v_0,v_1$ be the vertices of $\tau_0$, and suppose that we can
find, for each $i=1,\ldots,n$, an index $j(i)\in\{0,1\}$, such
that
$$\bar{\pi}_i(v_{j(i)})\in \bar{\pi}_i(v_{1-j(i)})+\R_{+}^n,$$
where $\bar{\pi}_i$ is the projection
$$\bar{\pi}_i:\oplus_{k=1}^{n}\R e_k\rightarrow \oplus_{k\neq i}\R e_k.$$
Then there would be two indices $i_1,i_2$ such that
$j(i_1)=j(i_2)$, wich would imply that already $v_{j(i_1)}\in
v_{1-j(i_1)}+\R_{+}^n$.
 This is impossible. Hence, we may
suppose that $\pi_2(v_j)\notin \pi_2(v_{1-j})+\R_{+}^{n-1}$,
 for $j=0,1$. By induction, we may assume that
$\pi_{n-1}(v_j)\notin \pi_{n-1}(v_{1-j})+\R_{+}^{2}$,
 for $j=0,1$.

This implies that $g(z_{n-1},z_n):=f_{\tau_0}(1,\ldots,1,z_{n-1},z_n)$
 satisfies
$g=g_{\tau_0(g)}$.
By Theorem \ref{assum}, the only thing left to prove, is that
we can find an index $i\in\{n-1,n\}$, such that
$\tilde{\tau}_0(f)+ \R e_i=\R^{n}$.
But if this were not the case, $\tilde{\tau}_0(f)$
would contain $\R e_{n-1}+\R e_n$,
and this contradicts that
 $\pi_{n-1}(v_j)\notin \pi_{n-1}(v_{1-j})+\R_{+}^2$, for $j=0,1$.
\end{proof}

To conclude, we give a proof of the Conjecture for $n=3$.
\begin{prop}
Let $f$ be a complex analytic singular germ at $0\in \C^3$,
let $\varphi$ be a Schwarz function on $\C^3$ with sufficiently small support,
 nonzero at $0$. We denote by $Z(s)$ the complex local zeta function
associated to $f$ and $\varphi$, and we suppose that $s_0\notin \Z$.
The polar multiplicity of $Z(s)$ at $s_{0}$ equals
$\rho$ if and only if $\tau_{0}$ is a stable face.
\end{prop}
\begin{proof}
The case $\rho=3$ is trivial, so we may assume $\rho<3$. By
Theorem \ref{cconj}, Proposition \ref{complexconj}, and Corollary
\ref{cor3}, we may furthermore suppose that $\tau_0$ is stable,
and not compact. As is easily seen, this implies that $t_0\geq 1$.
\end{proof}

\bibliographystyle{plain}
\bibliography{wanbib,wanbib2}
\end{document}